\documentstyle[12pt]{article}
\input mssymb
\title{Deformations of quadratic algebras and
the corresponding quantum semigroups}
\author{J. Donin and S. Shnider\\
{\normalsize Department of Mathematics, Bar-Ilan University}}

\begin{document}
\maketitle
\baselineskip=20pt
\parskip=5pt
\parindent=18pt
\setcounter{page}{1}
\def\hsp{\hspace*{18pt}}
\newcommand{\QED}{\hspace{0.2in}\vrule width 6pt height 6pt depth 0pt
\vspace{0.1in}}
\newcommand{\Cal}{\cal}
\newcommand{\Proof}{{\em Proof} \hspace{0.2in}}
\newcommand{\tag}[1]{\eqno(#1)}
\newcommand{\BC}{$\Bbb C$}
\newcommand{\BCn}{$\Bbb C^n$}
\newcommand{\CP}{$\Cal P$}
\newcommand{\CO}{$\Cal O$}
\newcommand{\HK}{$\hat K$}
\newtheorem{theorem}{Theorem}[section]
\newtheorem{lemma}[theorem]{Lemma}
\newtheorem{proposition}[theorem]{Proposition}
\newtheorem{corollary}[theorem]{Corollary}
\newtheorem{problem}[theorem]{Problem}
\newtheorem{definition}[theorem]{Definition}

\def\Un{{1\!\!{\rm l}}}
\def\un{{\rm 1\mkern-5mu  l }} 
\def\Rb{{I\!\! R}}
\def\Cb{\ \hbox{\vrule width 0.6pt height 6pt depth 0pt
                      \hskip -3.5 pt} C}

\def\t{\otimes}
\def\ff{\varphi}
\def\cI{{\cal I}}
\def\th{{\rm th}}
\def\k{{\bf k}}
\def\kh{\k[[h]]}
\def\la{\lambda}
\def\Hom{{\rm End}}
\def\End{{\rm End}}
\def\Aut{{\rm Aut}}
\def\Ad{{\rm Ad}}
\def\ff{\varphi}
\def\R{{\Cal R}}
\def\g{{\Cal G}}
\def\F{{\Cal F}}
\def\tb{{\bf t}}
\def\C{{\Cal C}}
\def\S{{\Cal S}}
\def\De{{\Delta}}
\def\ti{\tilde}
\def\I{{\Cal I}}
\def\K{{\Cal K}}

\begin{abstract}
Let $V$ be a finite dimensional vector space. Given a decomposition
$V\otimes V=\oplus_i^n I_i$, define $n$ quadratic algebras $(V, J_m)$
where $J_m=\oplus_{i\neq m} I_i$. This decomposition defines also
the quantum semigroup $M(V;I_1,...,I_n)$ which acts on all these quadratic
algebras. With the decomposition we associate a family of associative
algebras $A_k=A_k(I_1,...I_n)$, $k\geq 2$. In the classical case, when
$V\otimes V$ decomposes into the symmetric and skewsymmetric tensors,
$A_k$ coincides with the group algebra of the symmetric group $S_k$.
Let $I_{ih}$ be deformations of the subspaces $I_i$. In the
paper we give a criteria for flatness of the corresponding deformations
of the quadratic algebras $(V[[h]],J_{ih}$ and the quantum semigroup
$M(V[[h]];I_{1h},...,I_{nh})$. It says that the deformations will
be flat if the algebras $A_k(I_1,...,I_n)$ are semisimple and under
the deformation their dimension does not change.

Usually, the decomposition into $I_i$ is defined by a given Yang-Baxter
operator $S$ on $V\otimes V$, for which $I_i$ are its eigensubspaces,
and the deformations $I_{ih}$ are defined by a deformation $S_h$ of $S$.
We consider the cases when $S_h$ is a deformation of Hecke or Birman-Wenzl
symmetry, and also the case when $S_h$ is the Yang-Baxter operator
which appears by a representation of the Drinfeld-Jimbo quantum group.
Applying the flatness criteria we prove that in all these cases we
obtain flat deformations of the quadratic algebras and the corresponding
quantum semigroups.
\end{abstract}

\section{Quadratic algebras, quantum semigroup, and notations}
\label{s1}
Let $V$ be a module over a ring $A$, $I$ a submodule of
$V^{\t 2}=V\t V$.
Denote by $I^{i,k}$ the submodule in $V^{\t n}$ of the type
$V\t\cdots\t I\t\cdots\t V$ where $I$ occupies the
positions $i,i+1$. Set $I^k=\sum_i I^{i,k}$ and $I^{(k)}=\cap_i I^{i,k}$.
So, $I=I^2=I^{(2)}$.

Let $V^*$ be the dual module to $V$. Denote by $I^\bot$ the submodule
of $V^*$ which consists of all linear mappings $\ff:V\to A$ such that
$\ff(v)=0$ for $v\in I$.

We say that the ordered pair of submodules $(I,J)$,
$I,J\subset V^{\t 2}$, is well situated if
$V^{\t k}=I^{(k)}\oplus J^k$ for all $k\geq 2$. In particular,
$V^{\t 2}=I\oplus J$. It is easy to see that if the pair $(I,J)$
is well situated then the pair $(J^\bot, I^\bot)$ of submodules in
$V^{*\t 2}$ is also well situated. This follows from the relations
$(L+M)^\bot=L^\bot\cap M^\bot$ and $(L\cap M)^\bot=L^\bot+M^\bot$
which are true for any submodules of an $A$-module.

For a submodule $J\in V^{\t 2}$ we denote by $Q_J=(V,J)$ the quadratic
algebra $T(V)/\cI_J$, where $T(V)$ is the tensor algebra over $V$
and $\cI_J$ denotes the ideal generated by $J$. The algebra $Q_J$ is
a graded one, its $k^\th$ homogeneous component $Q_J^k$ is equal to
$T^k(V)/J^k$ as $A$-module. If the pair $(I,J)$ is well situated the
restriction of the natural mapping $T^k(V)\to Q_J^k$ gives an
isomorphism $I^{(k)}\to Q_J^k$ of $A$-modules.

In the sequel we will deal with the cases when $A$ is either a field $\k$ of
characteristic zero or the algebra $\k[[h]]$ of formal power series
in a variable $h$. In the latter case we will consider only modules of
finite rank and complete in $h$-adic topology. In particular,
all tensor products will be completed in that topology.
Any free $\kh$-module of rank $n$ is isomorphic as $\k[[h]]$-module
to $E\t_\k \kh=E[[h]]$, the module of formal power series in $h$
with coefficients from $E$.

We say that
a submodule $J_h$ of a $\kh$-module $E_h$ is a splitting submodule if it has
a complementary submodule $I_h$, i.e.
$E_h=J_h\oplus I_h$. It is clear that in case $E_h$ is a free module
any submodule $J_h$ is free, but
$J_h$ is a splitting one
if and only if the module $E_h/J_h$ is free.
We call a morphism of free modules, $\ff:E_h\to V_h$, flat
if $Im\ff$ (or equivalently, $Ker\ff$) is a splitting submodule.

Let $J$ is a linear subspace in a vector space $E$ over $\k$.
We say that $J_h$ is a family of subspaces in $E$, or a (formal)
deformation of the subspace $J$, if $J_h$ is a splitting submodule
in $E_h=E[[h]]$ such that $J_0=J$. Here $J_0$ is the set of elements
which obtained from elements of $J_h$ replacing $h$ by $0$.
Note that for a submodule $J_h$ in $E[[h]]$ we have
$\dim J_0\leq\dim J_{h\neq 0}$, and $J_h$
defines a deformation of the subspace $J_0$
if $\dim J_0=\dim J_{h\neq 0}$. Here $J_{h\neq 0}$ denotes the vector
subspace in the ``general'' point, i.e. the vector subspace
$J_h\t_{\kh}\k\{\!\{h\}\!\}$ in the vector space $V_h\t_{\kh}\k\{\!\{h\}\!\}$
over the field of formal Laurent series $\k\{\!\{h\}\!\}$.

If $J_h$ is not a splitting submodule in $E[[h]]$. then the module
$P_h=E[[h]]/J_h$ has a decomposition $P_h=P^{'}_h\oplus P^{''}_h$
where $P^{'}_h$
is a free module and $P^{''}_h$ is the torsion submodule of $P_h$,
that is $b\in P^{''}_h$ if and only if there exists $m>0$ such that
$h^m b=0$. Denote by $J^{'}_h$ the kernel of the natural projection
$E[[h]]\to P^{'}_h$. It is clear that $J^{'}_h$ is a splitting
submodule in $E[[h]]$ and
$\dim J_{h\neq 0}=\dim J^{'}_{h \neq 0}=\dim J^{'}_0$. So we denote
$J^{'}_0=J_{h\to 0}$.

We will only consider quadratic $\kh$-algebras
$(V_h,J_h)$ such that $V_h$ is a free module of finite rank and
$J_h$ is a splitting submodule of $V_h^{\t 2}$. We associate to
the quadratic algebra $(V_h,J_h)$ the quadratic algebra $(V,J)$
over $\k$ taking $V=V/hV_h$ and $J=J_h/hJ_h$ with the natural
imbedding $J\to V^{\t 2}$.
In this case we call
the quadratic algebra $Q_h=(V_h,J_h)$ a deformation of the algebra
$Q=(V,J)$. We call another deformation $(V^{'}_h,J^{'}_h)$
of the algebra $(V,J)$ equivalent to $(V_h,J_h)$ if there exists
an isomorphism $\phi:V^{'}_h\to V_h$ which induces the identity
isomorphism on $V$ and the isomorphism
$\phi\t\phi:V^{'}_h\t V^{'}_h\to V_h\t V_h$ gives an isomorphism
$J^{'}_h\to J_h$. Since there exists an isomorphism $V[[h]]\to V_h$,
any deformation $(V_h,J_h)$ of the algebra $(V,J)$ can be given
by a formal deformation of the subspace $J$ in $V^{\t 2}$, i.e.
is equivalent to a deformation of the form $(V[[h]],J^{'}_h)$
where $J^{'}_0=J$.

Let $(V_h,J_h)$ be a quadratic algebra. Note that in general
$J_h^k$ need not be a splitting
submodule in $V_h^{\t k}$ for $k>2$, so the homogeneous component
$Q_h^k=V_h^{\t k}/J_h$ will not be a free module. We call the
deformation $Q_h$ a flat deformation (or quantization) of $Q$  if
all modules $Q_h^k$ are free. Note that this terminology is not
completely standard. For many authors the flatness condition is
included in the definition of a deformation.

We mention here a theorem due to Drinfeld
\cite{Dr},
which states that in case of the Koszul quadratic algebra $(V,J)$,
\cite{Ma}, in order
for all modules $Q_h^k$, $k>2$, to be free it is sufficient that the module
$Q_h^3$ be free, i.e. the submodule $J_h^3$ be splitting in
$V_h^{\t 3}$.

Let $V$ be a vector space over $\k$. Suppose $I_i$, $i=1,...,n$, are
vector subspaces in $V^{\t 2}$ such that $V^{\t 2}=\oplus_i^n I_i$.
Denote $J_k=\oplus_{i\neq k}I_i$, so $V^{\t 2}=I_i\oplus J_i$.
We associate to the tuple $(I_1,...,I_n)$ a quantum semigroup
$M(V)=M(V;I_1,...,I_n)$ in the
following way. We identify $\End(V)=V\t V^*$ and put
$M(V)=(\End(V),I)$, the quadratic algebra where the subspace
$I$ of $\End(V)^{\t 2}$ is defined as
$I=\sigma_{2,3}(I_1\t I_1^\bot+\cdots +I_n\t I_n^\bot)$
where $\sigma_{2,3}$
is the permutation of the second and third tensor components.
The algebra $M(V)$ has the natural bialgebra structure and
the algebras $(V,J_i)$ make into comodules over $M(V)$~\cite{Ma}.

The quantum semigroup $M(V)$ also admits the description in the spirit
of Faddeev-Reshetikhin-Takhtajan \cite{FRT}. Let $\la_i$, $i=1,...,n$,
be different elements from $\k$. Let $S$ be the linear operator acting
on $\End(V^{\t 2})$, which has $I_i$ as the eigensubspace
corresponding to the eigenvalue $\la_i$ for all $i$.
Identifying $\End(V^{\t 2})\cong \End(V)^{\t 2}$ via the
Kronecker product we may view $S$ as an element of
$\End(V)^{\t 2}$, $S=S_{(1)}\t S_{(2)}$ in the Sweedler notation.
Then $I$ consists of the elements $X=X_{(1)}\t X_{(2)}$ of
$\End(V)^{\t 2}$ having the form
$SX-XS=S_{(1)}X_{(1)}\t S_{(2)}X_{(2)}-X_{(1)}S_{(1)}\t
X_{(2)}S_{(2)}$.
This means that coaction of $M(V)$ on $V$ preserves
all the subspaces $I_i$.

Denote by $A_2(S)$ the associative
subalgebra in $\End(V)$ generated by $S$. It is a
semisimple algebra isomorphic to a direct sum of $n$ copies of the
base field. Let $A_k(S)$ be the associative subalgebra in
$\End(V^{\t k})$ generated by the operators $S_i$, $i=1,...,k-1$,
where $S_i$ denotes the operator in $\End(V^{\t k})$ which coincides
with $S$ in the position $i,i+1$ and is the identity in the other
positions. It is clear that all the algebras $A_k(S)$
depend only on the subspaces $I_i$ but not on choosing of $\la_i$.
So the algebras $A_k(S)$ we also denote by $A_k(I_1,...,I_n)$.

\section{Quadratic algebras and semisimplicity}
\label{s2}
 In what following we suppose that the field $\k$ is equal to
$\Rb$ or $\Cb$.

A finite-dimensional representation  $E$ of an algebra $A$
(or $A$-module) is called simple if there are no nontrivial invariant
subspaces, and it is called
semisimple if $E$ is isomorphic to a direct sum of simple representations.
A finite-dimensional algebra is called semisimple if all its
finite-dimensional representations are semisimple.
An linear operator $B\in \End(E)$ is semisimple if the subalgebra
of operators generated by it is semisimple.
In general, we call a set of operators
${\cal F}\subset\End(E)$ semisimple if
the subalgebra $A({\cal F})$ generated by this family is semisimple.

As is known \cite{Pie} an algebra $A$ will be semisimple if and only
if its semisimple representations separate points, i.e. for any two elements
$a,b\in A$ there exists a semisimple representation $\varphi:A\to
\End(V)$ such that $\varphi(a)\neq\varphi(b)$. In particular,
if $A$ is a subalgebra of $\End(E)$ and the space $E$ is a semisimple
$A$-module then $A$ is semisimple. It follows from this that the
following algebras are semisimple:\\
\hsp a) $A(\varphi({\cal G}))$ for a representation
$\varphi:{\cal G}\to \End(E)$ of
a semisimple or compact Lie algebra ${\cal G}$;\\
\hsp b) $A(\varphi(G))$ for a representation $\varphi:{ G}\to \End(E)$ of
semisimple or compact Lie group ${ G}$;\\
\hsp c) $A({\varphi(\cal F}))$ for any subset ${\cal F}$ of a compact
Lie algebra or group and $\varphi$ is its representation.\\
Note that if $\varphi$ is a representation of a connected Lie group
$G$ and $\psi$ is the corresponding representation of its Lie algebra
${\cal G}$ then the algebras $A(\varphi(G))$ and $A(\psi({\cal G}))$
coincide.

\begin{proposition}
\label{p2.1}
Let $E$ be a finite-dimensional vector space over $\k$. Suppose $B_1,...,B_m$
is a semisimple set of linear operators on $E$ and $\la_1,...,\la_m$
are elements from $\k$. Denote
$L=\sum_iIm(B_i-\la_i)$, $K=\cap_iKer(B_i-\la_i)$. Then\\
\hsp a) the subspaces $L$ and $K$ are invariant under all the $B_i$;\\
\hsp b) $E=L\oplus K$.
\end{proposition}
\Proof
a) The invariance of $K$ is obvious. Let $v=\sum_i(B_i-\la_i)u_i$.
Then
$B_jv=\sum_iB_j(B_i-\la_i)u_i=\sum_i(B_i-\la_i)B_ju_i+
\sum_i[B_j,(B_i-\la_i)]u_i$, and a) follows from the equality
of commutators: $[B_j,(B_i-\la_i)]=[(B_j-\la_j),(B_i-\la_i)]$.\\
\hsp b) Because of semisimplicity there exists an invariant
subspace $P$ in $E$ complementary
to $L$. If $v\in P$ then $(B_i-\la_i)v$ has to belong to both
$L$ and $P$. Hence $(B_i-\la_i)v=0$ for all $i$.
It means that $v\in K$. So $P\subset K$ and $E=L+K$.
Let now $T$ be the invariant subspace in $E$ complementary to $K$.
It is clear that $L=\sum_i(B_i-\la_i)T$, so $L\subset T$
and, therefore, $L\cap K=0$. It implies that $E=L\oplus K$.~\QED

Now we consider deformations of semisimple algebras and their morphisms.
In general, let $A_h$ be an algebra over $\k[[h]]$ which is a free
$\k[[h]]$-module. Then $A_0=A_h/hA_h$ is an algebra over $\k$, and we call
$A_h$ a family of algebras, or a deformation of the algebra $A_0$.
If $A^{'}_h$ is another deformation of $A_0$ then a morphism of
the deformations is a $\k[[h]]$-algebra morphism $A_h\to A^{'}_h$
which is the identity for $h=0$.
The deformation is trivial if there exists an $\k[[h]]$-algebra
isomorphism $A_h\to A_0[[h]]=A_0\t_\k \k[[h]]$.
We say that a subalgebra $B_h\subset A_h$ is splitting if it is a
splitting $\k[[h]]$-submodule in $A_h$.

\begin{proposition}
\label{p2.2}
a) Let $A_h$ be a family of algebras.
Suppose the algebra $A_0$ over $\k$ is semisimple.
Then $A_h$ is isomorphic to $A_0[[h]]$ as $\k[[h]]$-algebra,
i.e. the deformation is trivial.\\
\hsp b) Let $\phi_h:A[[h]]\to B[[h]]$ be a morphism of $\k[[h]]$-algebras.
It induces the morphism $\phi_0:A\to B$ of $\k$-algebras.
Suppose $A$ is semisimple and $B$ is an arbitrary unital algebra.
Then there exists an element $f_h\in B[[h]]$ such that $f_0=1$
and $\phi_h=\Ad(f_h)(\phi_0\t \un)$. Here $\Ad(b)c=bcb^{-1}$ and
$\phi_0\t \un:A\t_\k \k[[h]]\to B\t_\k \k[[h]]$ is the morphism of
tensor products induced by $\phi_0$ and the identity morphism.
\end{proposition}
\Proof The proposition follows from the fact that the Hochschild
cohomology of any
semisimple algebra are equal to zero \cite{Pie} using the standard arguments
\cite{GGS}. More precisely, a) follows from $H^2(A,A)=0$ and b) from
$H^1(A,B)=0$ where $B$ is considered as $A$-bimodule via the morphism
$\phi_0$.~\QED

Let families of algebras $A_h$ and vector spaces $V_h$ be given.
Suppose the algebra $A_h$ acts on $V_h$, i.e. we are given a morphism
of $\k[[h]]$-algebras $\ff_h:A_h\to \End(V_h)$. Then $\ff_h$ induces
a morphism $\ff_0:A_0\to\End(V_0)$. On the other hand, any morphism
$\psi:A_0\to \End(V_0)$ generates in the trivial way the morphism
$\psi:A_0[[h]]\to \End_{\k[[h]]}(V_0[[h]])=\End_{\k}(V_0)[[h]]$,
so as a consequence of the preceding
proposition we get
that if the algebra $A_0$ is semisimple then the morphisms
$\ff_h$ and $\ff_0$ are isomorphic.

\begin{proposition}
\label{p2.3}
Let $E$ be a finite-dimensional vector space over $\k$. Suppose $B_1,...,B_m$
is a semisimple set of semisimple linear operators on $E$, $\la_1,...,\la_m$
are elements from $\k$. Let $B_{ih}\in \End(E)[[h]]$ and
$\la_{ih}\in \k[[h]]$  be deformations of
$B_i$ and $\la_i$ ($i=1,...,m$) such that\\
\hsp i) all the subalgebras $A_{ih}=A(B_{ih})$ are splitting submodules;\\
\hsp ii) the subalgebra $A_h=A(B_{1h},...,B_{mh})$ is a splitting submodule;\\
\hsp iii) all the submodules $K_{ih}=Ker(B_{ih}-\la_{ih})$ are
splitting ones.\\
Denote $L_h=\sum_iIm(B_{ih}-\la_{ih})$,
$K_h=\cap_iKer(B_{ih}-\la_{ih})$.
Then\\
\hsp a) the submodules $L_h$ and $K_h$ are ivariant under all
$B_{ih}$;\\
\hsp b) $E[[h]]=L_h\oplus K_h$. In particular, $L_h$ and $K_h$ are
splitting submodules.
\end{proposition}
\Proof The invariance of $L_h$ and $K_h$ can be proven as in
Proposition \ref{p2.1}.
At first, suppose that the algebra $A_h$ has the form $A_0[[h]]$ there
$A_0=A(B_1,...,B_m)$. Denote $L_0=\sum_iIm(B_{i}-\la_{i})$,
$K_0=\cap_iKer(B_{i}-\la_{i})$. If $K_0=0$ then $L_0=E$ by
Proposition \ref{p2.1}, and b) is obvious from the fact that
$\dim L_0\leq \dim L_{h\neq 0}$.
Suppose $K_0\neq 0$. Then, since
$K_0$ is an eigensubspace for all elements from $A_0$, the elements
$\la_i$ define an algebra homomorphism $\chi_0:A_0\to \k$ by
$\chi_0(B_i)=\la_i$. In the same way the element $\la_{ih}$ and
submodule $K_{ih}$ define a morphism $\rho_{ih}:A_{ih}\to\k[[h]]$
for all $i=1,...,m$.
Consider the morphism
$\chi_h=\chi_0\t\un:A_0[[h]]\to\k[[h]]$. Since the algebras $A(B_i)$ are
semisimple and the restriction of $\chi_0$ onto $A(B_i)$ coincides
with $\rho_{i0}$ for all $i$, it follows from i) and
Proposition \ref{p2.2}~b) that
$\chi_h=\rho_{ih}$ on $A_{ih}$ for all $i$. This implies that
$$L_h=\sum_{B\in A_h}Im(B-\chi_h(B)),\
K_h=\bigcap_{B\in A_h}Ker(B-\chi_h(B)).$$
Taking into account that $\chi_h=\chi_0\t\un$ we get that
$L_h=L_0[[h]]$ and $K_h=K_0[[h]]$, which proves the proposition in
the case $A_h=A_0[[h]]$.

Suppose now that $A_h$ is arbitrary. Then, by ii) and Proposition \ref{p2.2}
there exists an element $f_h\in \End(V)[[h]]$ such that
$f_hA_hf_h^{-1}=A_0[[h]]$. Constructing the spaces
$L_h^{'}=\sum_iIm(B^{'}_{ih}-\la_{ih})$,
$K_h^{'}=\cap_iKer(B^{'}_{ih}-\la_{ih})$ for
$B^{'}_{ih}=f_hB_{ih}f_h^{-1}$ we obtain that the modules
$L_h=f_h^{-1}L_h^{'}f_h$ and $K_h=f_h^{-1}K_h^{'}f_h$
satisfy the proposition.~\QED

Let $B_1,...,B_m$ is a semisimple set of semisimple operators in
a vector space $E$.
We say that deformations of these operators, $B_{1h},...,B_{mh}$,
form a flat deformation of the set if the conditions i) and ii)
from Proposition \ref{p2.3} hold.

{\bf Remark 2.1.} It is clear that if a semisimple operator $B$ on $E$
and its flat deformation $B_h$ are given then for any eigenvalue $\la$
of $B$ its deformation $\la_h$ is uniquely defined. Furthermore,
$K_h=Ker(B_h-\la_h)$ and $L_h=Im(B_h-\la_h)$ form deformations
of the subspaces $K=Ker(B-\la)$ and $L=Im(B-\la)$ in $E$.
Indeed, $\la$ defines a character $\chi:A(B)\to \k$, $\chi(B)=\la$, which,
by Proposition \ref{p2.2}, has the unique extension
$\chi_h:A(B_h)\to\k[[h]]$. Then, $\la_h=\chi_h(B_h)$.
So, it follows from this
that if $B_i$, $\la_i$,
$i=1,...,m$, is a set of semisimple operators on $E$ with fixed
eigenvalues and $B_{ih}$ is a flat deformation of the set,
then  the deformations of the eigenvalues, $\la_{ih}$, exist and
are uniquely
defined such that the condition (iii) of Proposition \ref{p2.3}
is satisfied and, therefore, for these $\la_{ih}$ the proposition
holds.

Let $\hat{E}$ be a tensor space over $E$, i.e. a tensor product
of a number of copies of $E$ and $E^*$.
The Lie algebra $\End(E)$ acts on $\hat{E}$ in the usual way.
For example, if $B$ is a linear operator
in $\End(E)$ and $v\t u\in E\t E$ then by definition
$\hat{B}(v\t u)=Bv\t u+v\t Bu$,
if $w\in E^*$ then $\hat{B}(w(v))=-w(B(v))$ for all $v\in E$.
In particular, if we identify $\End(E)\cong E\t E^*$ and
$M\in \End(E)$ then $\hat{B}(M)=BM-MB$.

\begin{proposition}
\label{p2.4}
Let $B_1$,...,$B_m$ be a semisimple set of semisimple linear operators
on $E$. Suppose that their deformations $B_{1h}$,...,$B_{mh}$ form a
flat deformation of the set.
Then for any tensor space $\hat{E}$ over $E$ the deformations
$\hat{B}_{1h}$,...,$\hat{B}_{mh}$ form a flat deformation of the set
$\hat{B}_{1}$,...,$\hat{B}_{m}$.
\end{proposition}
\Proof From semisimplicity it follows that there exists an element
$f\in \End(E)[[h]]$ such that the algebra
$\Ad(f)A(B_{1h},...,B_{mh})$ is equal to the algebra
$A(B_{1},...,B_{m})[[h]]$.
The group $\Aut(E)$ acts on $\hat{E}$ in the usual way.
Let $\hat{f}$ be the image
of $f$ by the corresponding mapping $\Aut(E)[[h]]\to \Aut(\hat{E})[[h]]$.
It is clear that the algebra $\Ad(\hat{f})A(\hat{B}_{1h},...,\hat{B}_{mh})$
is equal to the algebra $A(\hat{B}_{1},...,\hat{B}_{m})[[h]]$.
The algebra $A(\hat{B}_{1},...,\hat{B}_{m})$ is generated by the image
of the Lie subalgebra $L=L(B_1,...,B_m)$ of $\End(E)$ spanned on $B_1,...,B_m$.
The action of $L$ on $E$ and, therefore, on $\hat{E}$ is semisimple.
Hence, the algebra $A(\hat{B}_{1},...,\hat{B}_{m})$ is semisimple, which
proves the proposition.~\QED

{\bf Remark 2.2.} Proposition \ref{p2.4} remains true if the operators
$B_1$,...,$B_m$ are invertible and the operators
$\hat{B}_{1h}$,...,$\hat{B}_{mh}$
are defined as the images of $B_{1h}$,...,$B_{mh}$ by the mapping
$\Aut(E)[[h]]\to \Aut(\hat{E})[[h]]$.
To prove this we replace in the
proof above the Lie subalgebra $L=L(B_1,...,B_m)$ by the Lie subgroup
generated by $B_1,...,B_m$.

Now we return to the setting of the end of Section \ref{s1}.

\begin{proposition}
\label{p2.5}
Let $V$ be a finite-dimensional vector space over $\k$,
$S$ a linear operator on $V^{\t 2}$, and $S_h$ a deformation of $S$.
Suppose $\la_1,...,\la_n$ are the eigenvalues of $R$ and $I_1,...I_n$ are
the corresponding eigensubspaces.
Suppose that the subalgebras $A_k(S)$ in $\End(V^{\t k})$
are semisimple and the subalgebras $A_k(S_h)$ in $\End(V^{\t k})[[h]]$
are splitting submodules for all $k\geq 2$.
Then deformations of the eigenvalues,
$\la_{ih}$, and eigensubspaces, $I_{ih}$, $i=1,...,n$,
are uniquely defined and\\
\hsp a) the pairs of submodules $(I_{mh},J_{mh})$, $m=1,...,n$, are well
situated, where $J_{mh}=\oplus_{i\neq m}I_{ih}$;\\
\hsp b) the quadratic algebras $(V[[h]],J_{mh})$ form flat deformations
of the quadratic algebras $(V,J_m)$ for all $m$;\\
\hsp c) the quantum semigroup $M(V[[h]];I_{1h},...,I_{nh})$ is a flat
deformation of the quantum semigroup $M(V;I_1,...,I_n)$.
\end{proposition}
\Proof
It is clear that $S_{1h},...,S_{(k-1)h}$ form a flat set of semisimple
operators in $\End(V^{\t k})$ for all $k$.
The deformations of the eigenvalues,
$\la_{ih}$, and eigensubspaces, $I_{ih}$, $i=1,...,n$,
are uniquely defined by Remark 2.1.
Noting that
$I_m^{(k)}=\cap_{i=1}^{k-1}Ker(S_i-\la_m)$ and
$J_m^{k}=\sum_{i=1}^{k-1}Im(S_i-\la_m)$ and applying Proposition
\ref{p2.3} we obtain a) and b). Condition c) follows from
Proposition \ref{p2.4}~\QED

{\bf Remarks 2.3} 1. Proposition \ref{p2.5} gives the following criteria of
flatness for deformations of quadratic algebras and the corresponding
quantum semigroups.
Let $V$ be a vector space over $\k$. Suppose that $I_i$, $i=1,...,n$, are
vector subspaces in $V^{\t 2}$ such that $V^{\t 2}=\oplus_i^n I_i$.
Denote $J_k=\oplus_{i\neq k}I_i$, so $V^{\t 2}=I_i\oplus J_i$.
Suppose deformations $I_{ih}$, $i=1,...,n$, of the subspaces are given
and the subalgebras $A_k(I_{1h},...,I_{nh})$ of $\End(V^{\t k})[[h]]$ are
splitting submodules for all $k$. Then all the deformations $(V[[h]],J_{mh})$
of the quadratic algebras $(V,J_m)$ are flat. Moreover, the
deformation $M(V[[h]];I_{1h},...,I_{nh})$ of the quantum semigroup
$M(V;I_1,...,I_n)$ is flat as well.

2. One can consider the case when the variable $h$ runs through
a complex or real analytic manifold $X$,
the subspaces $I_{ih}$ depend on $h$ analytically,
and one has the decomposition $V^{\t 2}=\oplus_i^n I_{ih}$ at any
point $h\in X$. Suppose $\dim A_k(I_{1h},...,I_{nh})$ does not depend
on $h$ (this condition replaces the condition of splitting of the
subalgebra in the formal case). Suppose $A_k(I_{1h},...,I_{nh})$ is
semisimple at one point of $X$. Then there exists an analytic subset
$Y_k\subset X$ such that for $h\in X\backslash Y_k$ all the
algebras $A_k(I_{1h},...,I_{nh})$ are semisimple and isomorphic to
each other (cf. Proposition \ref{p2.2} a) ).
Following the arguments of this section one can prove
that for $h\in X\backslash\cup Y_k$ all the quadratic algebras
$(V,J_{mh})$ have the same dimension of their homogeneous components.
The same is true for the corresponding quantum semigroups.

\section{Applications}
\label{s3}
1). Let $V$ be a finite-dimensional vector space over $\k$
($\k=\Rb$ or $\Cb$). Let $S$ be an invertible linear operator on $V\t V$
with two eigenvalues $\la$ and $\mu$ satisfying the braid relation
(or quantum Yang-Baxster equation)
$$S_1S_2S_1=S_2S_1S_2 \tag{1}$$
on $V^{\t 3}$.
In this case the subalgebras $A_k(S)\subset V^{\t k}$ are
images of the Hecke algebras. The Hecke algebra $H_k(\la,\mu)$ is
defined as the quotient algebra of the free algebra
$T(x_1,...,x_{k-1})$ of $k-1$ variables by the relations
$$x_ix_{i+1}x_i=x_{i+1}x_ix_{i+1}, \ \  x_ix_j=x_jx_i \mbox{ for }
\mid i-j\mid\geq 2, \tag{2}$$
$$(x_i-\la)(x_i-\mu)=0. \tag{3}$$
It is known, \cite{Co}, that for almost all pairs $(\la,\mu)$
(excepting an closed algebraic
subset)
this algebra is semisimple and isomorphic
to the group algebra, $H_k$, of
the symmetric group (the case $\la=1$, $\mu=-1$).
Moreover, in a neighborhood of each point $(\la_0,\mu_0)$
this isomorphism can be chosen analytically dependent on $\la$ and $\mu$.
We suppose that the eigenvalues, $\la$ and $\mu$, of $S$ correspond
to the semisimple Hecke algebra.
In this case $S$ is called a Hecke symmetry.
Gurevich \cite{Gu} considered the case in details.
Now we consider deformations of the Hecke symmetry.

Let $S_h$ be a deformation of the operator $S$ satisfying the
relations
$$S_{h1}S_{h2}S_{h1}=S_{h2}S_{h1}S_{h2}, \tag{4}$$
$$(S_h-\la_h)(S_h-\mu_h)=0, \tag{5}$$
where $\la_h$ and $\mu_h$ are deformations of $\la$ and $\mu$.
Let us prove that in this case the subalgebras $A_k(S_h)$ for all
$k\geq 2$, are splitting. Indeed,
due to relations (4) and (5) there exists an algebra homomorphism
$\phi_h:H_k[[h]]\to \End(V^{\t k})[[h]]$ such that
$Im(\phi_h)=A_k(S_h)$. Using Proposition \ref{p2.2} we conclude
that the algebra $A_k(S_h)$ is isomorphic to $A_k(S)[[h]]$ and,
therefore, splitting.

2). We obtain the same result if $S$ satisfies the Birman-Wenzl
relations:\\
\hsp a) the braid relation (1);\\
\hsp b) the cubic relation $(S-\la)(S-\mu)(S-\nu)=0$ for $\la, \mu,
\nu\neq 0$;\\
\hsp c) $P_1S_2P_1=aP_1$, where $P=(S-\la)(S-\mu)$ and $a$ is a constant;\\
\hsp d) $P_1P_2P_1=bP_1$, where $b$ is a constant.\\
It follows from b) that $S$ has three eigenvalues and eigensubspaces.

In this case the subalgebras $A_k(S)\subset V^{\t k}$ are
images of the Birman-Wenzl (BW) algebras $BW_k$~\cite{BW}.
The algebra $BW_k$ is
defined as the quotient algebra of the free algebra
$T(x_1,...,x_{k-1})$ of $k-1$ variables by the relations (2) and
\begin{eqnarray*}
(x_i-\la)(x_i-\mu)(x_i-\nu)=0,\\
p_ix_{i\pm 1}p_i=ap_i,\\
p_ip_{i\pm 1}p_i=bp_i,
\end{eqnarray*}
where $p_i=(x_i-\la)(x_i-\mu)$. One can show that the constants $a$ and
$b$ are uniquely defined, and $a=\la\mu(\la+\mu)$,
$b=(\la+\mu)^2\nu^2$.
Note that in \cite{BW} BW algebras are defined by eleven relations, see
\cite{Ke} where it is proven that the algebra $BW_k$ can be defined
as above.

It is known, \cite{BW}, that for almost all triples
$\la$, $\mu$, and $\nu$ this algebra is semisimple and analytically
depended on $\la$, $\mu$, $\nu$.
We suppose that $\la$, $\mu$ and $\nu$ form such a triple.
In this case $S$ satisfying the relations a)-d) is called
a Birman-Wenzl symmetry.
So, in the case of BW symmetry
algebras $A_k(S)$ are semisimple as well.

If by a deformation of $S$ the relations a)-d) hold
(with deformed eigenvalues $\la$, $\mu$, $\nu$) we say that this is
a deformation of Birman-Wenzl symmetry.
Using the same arguments as in 1) we obtain that by a deformation $S_h$ of
the BW symmetry $S$ the algebras $A_k(S_h)$ are splitting.

Applying Proposition \ref{p2.5} we obtain

\begin{proposition}
\label{p3.1}
Let $S$ be a Hecke (BW) symmetry on the space $V$, $I$ its
eigensubspace in $V\t V$, and $J$ is the sum of other eigensubspaces.
Suppose $S_h$ is a deformation of the
Hecke (BW) symmetry. Then the deformation defines a flat deformation
$(V[[h]],J_h)$ of the quadratc algebra
$(V,J)$ and the pair $(I_h,J_h)$ is well situated.
Moreover, the deformation of the quantum semigroup corresponding to
the eigensubspaces of $S$ is flat.
\end{proposition}

This proposition for the case $S=\sigma$ is proven in \cite{GGS1}.
Note that Gurevich proved in \cite{Gu} that in  case of Hecke symmetry
the algebra $(V,I)$ is Koszul. He also constructed Hecke symmetries
with nonclassical dimensions of homogeneous components of $(V,I)$.

In particular, deformations of the Hecke and BW symmetries
appears in \cite{FRT} by
construction of the quantum analogs (deformations) of the classical Lie groups.
Namely, the Hecke symmetry corresponds to the case of general linear
group, while the BW symmetry corresponds to the orthogonal and symplectic
cases.

3). Let $S$ be a Yang-Baxter (YB) operator on $V\t V$, i.e
$S$ is invertible and satisfies the braid relation (1).
Let $\hat{V}=U\t\cdots\t W$ be a tensor space over $V$, the spaces
$U$,...,$W$ are equal to $V$ or $V^*$. The group $\Aut(V)$ acts on
$\hat{V}$ in the usual way. Denote by $\hat{B}$
the image of $B\in \Aut(V)$ by the corresponding homomorphism
$\Aut(V)\to\Aut(\hat{V})$. The operator $S$ defines
the operator $\hat{S}=\hat{S}_{(1)}\t\hat{S}_{(2)}$ on
$\hat{V}\t\hat{V}$. Here we use the Sweedler notation,
$S=S_{(1)}\t S_{(2)}$. It is easy to see that $\hat{S}$ also
satisfies the braid relation.

Suppose now that $S$ is a Hecke symmetry. The operator $\hat{S}$
will not be a Hecke symmetry (it may have more than two eigenvalues),
but all the algebras $A_k(\hat{S})$ are semisimple.
Let $S_h$ be a deformation of $S$. This deformation defines a deformation,
$\hat{S}_h$, of $\hat{S}$. By 1) $\hat{S}_h$ defines flat deformations
of the algebras $A_k(S)$. So, using Remark 2.2 and Proposition
\ref{p2.5}, we obtain flat deformations of the quadratic algebras and
the quantum semigroup corresponding to the decomposition of
$\hat{V}\t\hat{V}$ into eigensubspaces of $\hat{S}$.
Of course, the similar statment is fulfilled for Birman-Wenzl
symmetry $S$.

4).
Let $U_h(\g)$ be the Drinfeld-Jimbo quantized universal enveloping
algebra (DJ quantum group), for a semisimple Lie algebra $\g$ over $\k=\Cb$.
Let $R\in U_h(\g)\t U_h(\g)$ be the corresponding quantum R-matrix.
Suppose, a representation $V_h$ of $U_h(\g)$ is given, which is a
deformation of the finite dimensional representation $V$ of $U(\g)$,
so $V_h$ is isomorphic to $V[[h]]$ as $\k[[h]]$-module, and the
representation can be presented as a homomorphism
$\rho:U_h(\g)\to \End(V)[[h]]$.
Consider the operator $S_h=\sigma R_h$ where
$R_h=(\rho\t\rho)(R)\subset\End(V)^{\t 2}[[h]]$ and $\sigma$ is the
standard permutation.
It is known that $S_h$ is a Yang-Baxter operator, i.e. satisfying the
braid relation (1). But it is not necessarily a flat deformation of
semisimple operator,
because at $h=0$ the operator $S_0$ is equal to $\sigma$,
so has two eigenvalues, $\pm 1$, while at the general point $h\neq 0$ it
is semisimple but may have more than two eigenvalues,
$\la_{ih}$, $i=1,...,n$, such that $\la_{i0}=\pm 1$. Nevertheless,
there is the decomposition $V[[h]]=\oplus_i I_{ih}$ where $I_{ih}$ are
eigensubmodules of $S_h$ corresponding to $\la_{ih}$, and, therefore,
all $I_{ih}$ are splitting. We will prove that also in this setting
the decomposition defines
flat deformations of quadratic algebras, $(V[[h]],J_{mh})$,
$J_{mh}=\oplus_{i\neq m}I_{ih}$, and of the corresponding quantum semigroup,
$M(V[[h]];I_{1h},...,I_{nh})$.
For this, according to Remark 2.3.1 we will show that the algebras
$A_k(I_1,...,I_n)$
are semisimple, $I_i=I_{i0}$, and $A_k(I_{1h},...,I_{nh})$ are splitting
for $k\geq 2$.

We recall some results of Drinfeld from \cite{Dr1} and \cite{Dr2}.
Additional structures on the category $Rep_A$ of representations of an
associative algebra $A$ and morphisms of these structures can be given by
the additional structures on the algebra $A$ itself. Thus, the structure of
quasitensor monoidal category on $Rep_A$ can be given with the help of
an algebra homomorphism $A\to A\otimes A$ (comultiplication), an element
$\Phi\in A^{\otimes 3}$ (associativity constraint), and R-matrix
$R\in A^{\otimes 2}$ (commutativity constraint),
satisfying the certain conditions. A morphism of such two
structures can be given by an element $F\in A^{\otimes 2}$. Drinfeld defined
such a structure on $A=U(\g)[[h]]$ for any semisimple Lie algebra $\g$ with the
usual comultiplication $\De$ but nontrivial $R$ and $\Phi$. He then proved
that the
corresponding quasitensor category is isomomorphic by some $F_h$ to the
category
of representations of the Drinfeld-Jimbo quantum group $U_h(\g)$
which coincides with
$U(\g)[[h]]$ as an algebra but has a noncommutative comultiplication
$\ti{\De}$. We denote the corresponding quasitensor categories by
${\C}$ and $\ti{\C}$, respectively. We keep the notations $\R_h$ and
$\Phi_h$ for the R-matrix and the associativity constraint in the
category $\C$, while $R_h$ denote the R-matrix for $\ti{\C}$.
Further, Drinfeld proved that $\R_h$ and $\Phi_h$ may be chosen as
$\R_h=e^{h\tb}$ where $\tb\in \g\t\g$
is the split Casimir, and $\Phi_h=e^{L(h\tb_1,h\tb_2)}\in U(\g)^{\t 3}[[h]]$
where $L(h\tb_1,h\tb_2)$ is a Lie expression of $\tb_1=\tb\t 1$
and $\tb_2=1\t\tb$. The element $F_h\in U(\g)^{\t 2}[[h]]$ is
congruent to $1\t 1$ modulo $h$ and satisfies the equation
$$(F_h\t 1)\cdot(\De\t id)(F_h)=(1\t F_h)\cdot(id\t\De)(F_h)\cdot\Phi_h.
\tag{6}$$
According to this, the commutativity constraints in the categories $\C$
and $\ti{\C}$ are given by the elements
$$\S_h=\sigma\R_h=\sigma e^{h\tb} \ \ \mbox{and} \ \
S_h=F_h\S_h(F_h)^{-1}=F_h\sigma e^{h\tb}(F_h)^{-1}, \tag{7}$$
respectively.

Let us come back to the setting from the beginning of 4).
Using (7) we obtain that the eigenvalues of the operators $\S_h$ and
$S_h$ acting on $(V\t V)[[h]]$ are $\la_{ih}=\pm e^{h\la_i}$ where
$\la_i$ are the eigenvalues of $\tb$ on the same space.

Now we transfer the setting to the category $\C$, i.e. we consider
$V[[h]]$ as an object of $\C$. Instead of $S_h$ we consider $\S_h$,
but the tensor products depend on the placement of parentheses and
the connection between two bracketing,
$\phi_h:(V^{\t k}){'}\to(V^{\t k}){''}$, of the $k$-fold tensor product is
generated by the operator
$$\Phi_h:((V\t V)\t V)[[h]]\to(V\t(V\t V))[[h]] \tag{8}$$
and looks like an expression $\phi_h\in \End(V^{\t k})[[h]]$
depending on the elements
$\tb_{i,j}=1\t\cdots\t\tb_{(1)}\t\cdots\t\tb_{(2)}\t\cdots\t 1$
($\tb_{(1)}$ and $\tb_{(2)}$ at the places $i$ and $j$,
$\tb=\tb_{(1)}\t\tb_{(2)}$ in the Sweedler notations).
It is easy to see that the eigensubmodules of
the operator $\S_h$ have the form $\I_i[[h]]$ where $\I_i$ are the
common eigensubmodules of $\tb$ and $\sigma$. Denote by
$A_k(\sigma,\tb)$ the subalgebra in $\End(V^{\t k})$ generated by the
elements $\sigma_{i,i+1}$ and $\tb_{i,i+1}$.
So we get that the
algebra $A_2(\I_1[[h]],...,\I_n[[h]])$ is equal to $A_2(\sigma,\tb)[[h]]$.
The algebras $A_k{'}(\I_1[[h]],...,\I_n[[h]])$ and
$A_k{''}(\I_1[[h]],...,\I_n[[h]])$ for two bracketings
$(V^{\t k}){'}$ and $(V^{\t k}){''}$ are connected by
$$A_k{'}(\I_1[[h]],...,\I_n[[h]])=\phi_h^{-1}A_k{''}(\I_1[[h]],...,\I_n[[h]])
\phi_h.$$
Using the relations of the type
$\tb_{1,3}=\sigma_{2,3}\tb_{1,2}\sigma_{2,3}$
we conclude that $\phi_h\in A_k(\sigma,\tb)[[h]]$,
and using the fact that $\phi_h$ is congruent to $1\t 1\t 1$ modulo $h^2$,
we conclude by induction on $k$ that $A_k{'}(\I_1[[h]],...,\I_n[[h]])=
A_k{''}(\I_1[[h]],...,\I_n[[h]])=A_k(\sigma,\tb)[[h]]$.

Passing to the category $\ti{\C}$ we obtain that
$A_k(I_{1h},...,I_{nh})=f_h{'}A_k{'}(\I_1[[h]],...,\I_n[[h]])(f_h{'})^{-1}
=f_h{'}A_k(\sigma,\tb)[[h]](f_h{'})^{-1}$ where $f_h{'}$ is the
composition (depending on the bracketing ${'}$)
of a number of $F_h$ with $\De$ applied appropriate factors.
We will obtain the same result applying to $A_k(\sigma,\tb)$ the
element $f_h{''}$ related to the bracketing ${''}$.
So, we have proved that $A_k(I_{1h},...,I_{nh})$ is splitting.

For $h=0$ this algebra is equal to $A_k(\sigma,\tb)$. Let us prove that
it is a semisimple algebra.
Indeed, $\tb$ may be presented as $\tb=\sum_id_i\t d_i$ where $d_i$
form an orthogonal (with respect to the Killing form) basis in the
maximal compact subalgebra $\K$ of $\g$. Hence, there exists a
Hermitian metric on $V$ invariant under action of $\K$. This metric
induces naturally the metric on $V^{\t k}$ which will be invariant
under the operators $\tb_{i,j}$ and $\sigma$. So these operators
are unitary ones, therefore the algebra $A_k(\sigma,\tb)$
generated by them is semisimple.

Applying  proposition \ref{p2.5} we obtain
\begin{proposition}
\label{p3.2}
Let $S_h$ be the Yang-Baxter operator on a space $V$,
which is obtained by the representation of the Drinfeld-Jimbo quantum group.
Then $S_h$ defines the decomposition $V^{\t 2}=\oplus_{i=1}^n I_{ih}$
into eigensubmodules. Let $J_{mh}=\oplus_{i\neq m} I_{ih}$.
Then
$(V[[h]],J_{mh})$ are flat deformations of the quadratc algebras
and the pairs $(I_{mh},J_{mh})$ are well situated.
Moreover, the deformation of the quantum semigroup
$M(V[[h]];I_{1h},...,I_{nh})$ is flat.
\end{proposition}

Another proof of this proposition is contained in \cite{DS}.
\bigskip
{\bf Acknowledgments.} The authors are grateful to J. Bernstein and
D. Gurevich for their interest to the paper and very helpful discussions.


\small
\begin{thebibliography}{GGS1}
\bibitem[BW]{BW} J.S  Birman and H. Wenzl, {\em Braids, Link
Polynomials and a New Algebra\/}, Trans. of AMS, v. 313, n. 1(1989),
pp. 249-273.

\bibitem[Co]{Co} M. Couillens, {\em Algebres de Hecke\/},
Seminaire sur les groupes finis, Tome II,
Publications math\'{e}matiques de
l'Universit\'{e} Paris VII.

\bibitem[Dr]{Dr} V.G. Drinfeld,
{\em On quadratic commutator relations in the quasiclassical case\/},
Math. physics and functional analysis, Kiev, Naukova Dumka, 1986.

\bibitem[Dr1]{Dr1} V.G. Drinfeld, {\em Quasi-Hopf algebra},
Leningrad Math. J., 1 (1990), pp.1419-1457.

\bibitem[Dr2]{Dr2} V.G. Drinfeld, {\em On Quasitriangular Quasi-Hopf
Algebras and a Sertain Group Closely Connected with
$Gal(\tilde{Q}/Q)$\/},
Leningrad Math. J., 2 (1990), pp.149-181.

\bibitem[DS]{DS} J. Donin and S. Shnider,
{\em Quasi-associativity and Flatness Criteria for Quadratic Algebra
Deformation\/}, Israel J. of Math., to appear.

\bibitem[FRT]{FRT} L.D. Faddeev, N.Yu. Reshetikhin, and L.A. Takhtajan,
{\em Quantization of Lie groups and Lie algebras},
Leningrad Math J., 1 (1990), pp. 193-225.

\bibitem[GGS]{GGS} M. Gerstenhaber, A. Giaquinto, and S.D. Shack,
{\em Algebraic cohomology and deformation theory\/},
Deformation theory of algebras and structures and applications
(M.Gerstenhaber and M. Hazewinkel, eds.),
Kluwer, Dordrecht, 1988, pp. 11-264.

\bibitem[GGS1]{GGS1} M. Gerstenhaber, A. Giaquinto, and S.D. Shack,
{\em Constructions of Quantum Groups from Belavin-Drinfeld
infinitesimals\/}, Quantum Deformations of Algebras and Their
Representations (A.Joseph and S.Shnider, eds.),
Israel Mathematical Conference Proceedings, 1993.

\bibitem[Gu]{Gu} D.I. Gurevich,
{\em Algebraic aspects of the quantum Yang-Baxter equation\/},
Leningrad Math. J., 2 (1991), pp. 801-828.

\bibitem[Ke]{Ke} S.V. Kerov, {\em Realization, Representations, and
Characters of the Birman-Wenzl algebra\/}, preprint.

\bibitem[Ma]{Ma} Yu. Manin,
``Quantum Groups and Non-commutative  Geometry'',
Les Publications Centre de R\'echerches Math\'ematiques,
Montr\'eal, 1989.

\bibitem[Pie]{Pie} R.S. Pierce, ``Associative Algebras'',
Springer-Verlag, New York Heidelberg Berlin, 1982.
\end{thebibliography}
\end{document}